\documentstyle[12pt,amscd]{amsart}

\newcommand{\Cal}{\cal}

\newtheorem{th}{Theorem}[section]

\theoremstyle{definition}

\newtheorem{df}{Definition}[section]

\theoremstyle{remark}

\numberwithin{equation}{section}

\newenvironment{ack}{\small \trivlist \item[\hskip \labelsep{\it
Acknowledgments}.]}{\endtrivlist}

\begin{document}

\title[Homological Reduction of Constrained Poisson Algebras]{Homological Reduction of Constrained Poisson Algebras}

\author {Jim Stasheff}  
\address
{Department of Mathematics, University of North Carolina, Chapel Hill, NC
27599-3250, USA}
\email{jds@@math.unc.edu}
\thanks {Research supported in part by NSF
grants DMS-8506637, DMS-9206929, DMS-9504871, a grant from the Institute for 
Advanced Study and a Research and Study Leave from the University of
North Carolina-Chapel Hill.
Announced in the Bulletin of the American Mathematical Society as 
      ``Constrained Poisson algebras and strong homotopy representations'' 
      [S2].  This paper includes the mathematical version of the physics 
in [FHST]}

%\date{January 14, 1996}   

\maketitle
%remove most of what follows upo to text

%\noindent
%{\bf {\Large Homological Reduction of Constrained Poisson Algebras}}
%\bigskip

\renewcommand{\thefootnote}{\arabic{footnote}}
\setcounter{footnote}{0}

%\hfill {\tt  }

%\maketitle

%\bigskip

%\NoBlackBoxes
\def\hatV{V/\cal F}
\def\c{\cite}
\def\smooth#1{C^\infty({#1})}                   % (real) smooth functions
\def\smoothc#1{C^\infty_{\complex}({#1})}       % (complex) smooth functions
\def\symmetric{\Cal{S}}                         % symmetric algebra
\def\tensor{\otimes}                            % tensor product
\def\iso{\approx}
\def\congruent{\equiv}
\def\alt{  $Alt_{P/I} (I/I^2,P /I)$ } 
\def\cite#1{[{\bf #1}]}
%\def\myspace{\hskip 1em\relax}
%\def\myhang{\hangvskip.5ex=4em\hangafter=0}
%\def\doit#1{\hskip-4em\hbox to 4em{#1\hss}}
%\def\refvskip.5ex#1{\par\noindent\doit{#1\myspace}\myhang}
%\hyphenation 
%
%\baselineskip=1.5\baselineskip
\hsize=6.5in 
\vsize=9in   
%\voffset=-.25in
%\pageno=1
%\headline={\ifnum\pageno>1 \tenbf Homological reduction \hss Stasheff\else\hfil\fi}
%

           Reduction of a Hamiltonian system with symmetry and/or 
constraints has a long history.  There are several reduction 
procedures, all of which agree in ``nice'' cases \c{AGJ}.  Some have 
a geometric emphasis - reducing a (symplectic) space of states
\c{MW}, while others are algebraic - reducing a (Poisson) algebra of 
observables \c{SW}.  Some start with a momentum map whose 
components are constraint functions \c{GIMMSY}; some start with a gauge 
(symmetry) algebra whose generators, regarded as vector fields, 
correspond via the symplectic structure to constraints \c{D}.  
      The relation between symmetry and constraints is particularly 
      tight in the case Dirac calls ``first class''.  The present paper 
      is concerned entirely with this first class case and deals with 
      the reduction of a Poisson algebra via homological methods, 
      although there is  considerable motivation from topology, 
      particularly via the models central to rational homotopy theory.  
           
Homological methods have become increasingly important in 
      mathematical physics, especially field theory, over the last 
      decade.  In regard to constrained Hamiltonians, they came into 
      focus with Henneaux's Report \c H on the work of Batalin, Fradkin 
      and Vilkovisky \c{BF,BV 1-3}, emphasizing the acyclicity of a 
      certain complex, later identified by Browning and McMullan as the Koszul 
      complex of a regular ideal of constraints.  I was able to put the 
      FBV construction into the context of homological perturbation 
      theory \c {S1} and, together with Henneaux et al \c {FHST}, extend the 
      construction to the case of non-regular geometric constraints of first 
      class.  Independently,  using a mixture of homological and 
      $C^1$ -patching techniques, Dubois-Violette extended the construction 
      to regular but not-necessarily-first-class constraints \c{D-V}.

I am grateful to all of the above for their input and inspiration,
whether in their papers or in conversation.  The present version
has also profitted from conversations at the MSRI Workshop on
Symplectic Topology.
Finally, I would like to express my thanks to the referee
who has read several versions with extreme care, suggesting extensive
improvements, both factual and stylistic.
While early revision was
in progress, Kimura sent me a copy of \cite {Ki} which has also had
a significant influence on the present exposition, as has his continued 
interaction while with me at UNC as an NSF Post-Doc.

      \section{   Preliminaries}
           This research touches on questions which it is hoped will be 
      of interest to mathematical physicists, symplectic and algebraic 
      geometers and homotopy theorists.  The techniques used here 
      are primarily those of differential commutative algebra and 
      rational homotopy theory.  We write with a dual vision and
hopefully a dual audience; for example, the constraints are functions
on a symplectic manifold and the physics literature speaks almost entirely
in terms of the constraints whereas the algebra can be expressed more
invariantly in terms of the ideal generated by the constraints.  
We work entirely over the reals $\Bbb R$ as our ground field, although
any field of characterisitic $0$ would do and the complex numbers $\Bbb C$
are more common in certain physical applications. The major 
Theorem 4.2 is expressed in algebraic terms, followed by remarks specifically 
in terms of the constraints themselves.

We begin therefore with a brief 
      (very!) review of the motivating background:  a tiny bit of 
      symplectic geometry, slightly more of Poisson algebra and the 
      essentials of constraint varieties and their symmetries in the 
      first class case.  The reader who desires more extensive 
      background or a more leisurely exposition may consult a variety  
      of sources listed in the bibliography.   The relations between the 
algebra and the motivating geometry are exposed particularly clearly in \c {Ki}.

      \subsection{The Hamiltonian Formalism}
           The motivating physical systems are described as 
      differential equations of motion or evolution involving smooth 
      functions on a manifold.  The underlying manifold W is assumed to 
      be symplectic.  This means there is a 2-form  $\omega$  such 
      that  $d\omega = 0$  and $\omega^{dim W} \neq 0$.  Equivalently, 
$\omega$  induces an isomorphism
      $$
TW  \rightarrow  T^*W.
$$
      (With an eye to future applications, we would like to allow  $W$  
      to be infinite dimensional, in which case the appropriate 
     definition is that the induced map  $TW \rightarrow T^* W$  be one-to-one.)
      In local coordinates  $q^1 ,...,q^n , p_1 ,...,p_n $, the form  
$\omega$  looks like $dq^i \wedge dp_i$ 
      (the summation convention will be assumed throughout this paper).
           
From an algebra point of view, the crucial point is 
      two-fold:  For any function  $f \in C^{\infty} (W)$, there is a 
Hamiltonian vector field  $X_f$   defined by  $\omega (X_f,\ ) = df$.
  For two functions  $f,g \in C^{\infty}(W)$, their Poisson bracket  $\{f,g\}
\in  C^{\infty} (W)$  is defined by
$$
                    \{f,g\} = \omega(X_f ,X_g ) = df(X_g ) = -dg(X_f ) .
 $$                          
      This bracket makes  $C^{\infty} (W)$  into a Poisson algebra, that is, a 
      commutative algebra  $P$  (with product denoted $fg$) together with a 
      bracket  $\{\  ,\  \} : P \otimes  P \rightarrow P$ forming a Lie algebra 
such that  $\{f,\ \}$ is a derivation of $P$ as a commutative algebra: 
 $\{f,gh\} = \{f,g\}h + g\{f,h\}$.  

A typical Hamiltonian system is one of the form $\{f,H\} =df/dt$ for fixed $H$.
Symmetries of such a system are given by functions $g$ which Poisson commute 
with $H$.  They form a sub-Lie algebra of $C^{\infty}(W)$.  
Symmetries arise also in connection 
with ``constraints''.  Regarded as in a dynamical system, solutions can be 
constrained to lie in a sub-manifold $V \subset  W$ (more generally, $V$ is 
just a sub-space), hereafter called the {\bf constraint locus}, also
known in the literature as a constraint surface.  As in algebraic geometry, 
we can think of $V$ as the zero 
set of some functions $\phi_\alpha :W \rightarrow R$, called {\bf constraints}. 
The algebra of functions $C^{\infty}$-in-the-sense-of-Whitney on $V$ can be 
identified with $C^{\infty} (W)/I$ where $I$ is the ideal of functions which  
vanish on $V$.  If $V\subset W$ is a closed and embedded submanifold, this
agrees with the usual notion of smooth functions on $V$.  
%Since for us $V$ will occur only as a subset of $W$, 
%we will abuse the language and refer to $C^{\infty} (W)/I$ as $C^{\infty}(V).$ 
%Henceforth, ``in-the-sense-of-Whitney'' will be understood where necessary.
	
Now if $W$ is symplectic (or just given a Poisson bracket on 
$C^{\infty} (W))$,
Dirac calls the constraints {\bf first class} if $I$ is closed under the Poisson
bracket.  (If the  ${\Bbb R}$-linear span of the $\phi_{\alpha}$   
is closed under the bracket, physicists say the  $\phi_{\alpha}$   
{\bf close} on a Lie algebra; this is a very nice case, but the more general 
      first class case is where homological techniques are really 
      important.)  When  the constraints are first class, we have that 
      the Hamiltonian vector fields $X_{\phi_\alpha}$  
  determined by the constraints are tangent to $V$ (where $V$ is smooth) and 
give a foliation $\cal F$ of $V$.  Similarly, $C^{\infty}(W)/I$ is a 
Lie module over $I$ with respect to the Poisson bracket.
In symplectic geometry, when $V$ is smooth, it is usually called a 
{\bf coisotropic} submanifold (see \c{W} for generalizations when $V$ is not 
smooth).  For the general case, we will call the constraint locus 
{\bf coisotropic} if the ideal is first class.

In many cases of interest, $I$ does not arise from the Lie 
      algebra of some Lie group of transformations of  $W$  or even  $V$, 
         but the corresponding Hamiltonian vector fields  $X_{\phi_\alpha}$  
are still referred to as (infinitesimal) symmetries.
In the nicest case, e.g. when the foliation $\cal F$  is given by a 
principal G-bundle structure on a smooth $V$, the algebra $C^{\infty}(V/\cal F)$
can be identified with the $I$-invariant sub-algebra of $C^{\infty}(W)/I$.  In 
great (if not complete) generality, this $I$-invariant sub-algebra represents 
the true observables of the constrained system.  

In this context, the ``classical BRST construction'', at least as developed by 
Batalin-Fradkin-Vilkovisky and phrased in terms of constraints, is a 
homological construction for performing
the reduction of the Poisson algebra $\smooth W$ of smooth functions on a
Poisson manifold $W$ by the ideal $I$ of functions which vanish on a
coisotropic constraint locus.  But the construction produces cohomology in other
degrees than zero, which at least in some cases, admits a geometric 
interpretation.

	Instead of considering just the ``observable'' functions, one can 
consider the deRham complex of longitudinal or vertical forms of the foliation
$\cal F$, that is, the complex $\Omega (V,\cal F)$ consisting of forms on 
vertical vector fields, those tangent to the leaves.
If we think of $\cal F$ as an involutive sub-bundle
of the tangent bundle to $V$, then $\Omega (V,\cal F)$ 
consists of sections of $\Lambda^*\cal F$.
  In adapted local coordinates $(x^1,...,x^{r+s})$ 
with $(x^1,...,x^r)$ being coordinates on a leaf, a typical longitudinal 
form is
$$
f_J (x) dx^J \quad \text { where} \quad J = (j_1,...,j_q) \quad \text{ with} 
\quad 1 \leq j_1 < ... j_q \leq r, \text{ \ the leaf dimension}.
$$
The usual exterior derivative of differential forms restricts to determine
the vertical exterior derivative because $\cal F$ is involutive.
This complex is familiar in foliation theory, c.f. \c{HH}.  The classical
BRST-BFV construction has, in the nice cases, the same cohomology as this
complex of longitudinal forms.

A  major motivating example for the BFV construction was provided by 
gauge theory.  Here $W$ is $T^*\cal A$ where $\cal A$ is the space of connections for a fixed principal $G$-bundle $G \to P \to B$.  The reduced phase space is $T^*(\cal A/\cal G)$ where $\cal G$ is the group of ``gauge transformations'',
the vertical automorphisms of $P$.

In considering what the physicists \c{BF},\c{BV1-3},\c {FF}, \c{FV},\c{H},\c{BM}
did in some special cases, I recognized a homological  
``model'' for    $\Omega(V, \cal F)$ in roughly the sense of rational
	homotopy theory  \c{Su}.  This is the same sense in which
      the Cartan-Chevalley-Eilenberg complex \c{CE} for the cohomology of a Lie 
      algebra  $\frak g$  is a ``model'' for  $\Omega ^* (G)$  where  $G$  
is a compact Lie group with Lie algebra  $\frak g$.  The physicists'
model is itself 
crucially a Poisson algebra extension of a Poisson algebra  $P$  and its 
      differential contains a piece which reinvented the Koszul 
      complex for the ideal  $I$.  The differential also contains a 
      piece which looks like the Cartan-Chevalley-Eilenberg 
      differential. Generalizations of the Cartan-Chevalley-Eilenberg 
differential as they 
occur in physics are usually referred to as BRST operators.  This honors 
seminal work of Becchi, Rouet and Stora \c {BRS} and, independently, 
Tyutin \c{Ty}.  Apparently it was the search for such an operator in aid of 
quantization which motivated the work of Batalin, Fradkin and Vilkovisky.

It was Browning and McMullan \c{BM} who first identified the Koszul complex 
      within the construction in the regular case, (Henneaux had already 
	called attention to the relevance of that acyclicity) leading both 
      Dubois-Violette \c{D-V} and myself \c{S1} independently to adopt a 
      more fully homological approach, although with somewhat different 
      emphases.  Dubois-Violette retains some of the symplectic 
      geometry and is able to handle regular general (not necessarily 
      first class) constraints.  On the other hand, by restricting to 
      first class constraints, 
	in joint work with Henneaux et al \c{FHST}, I was able to handle 
	non-regular ideals in suitable geometric circumstances.  

In the present paper, I start at the level of
the purely (Poisson) algebraic structures. In particular, I 
      adapt the notion of ``model'' from rational homotopy theory and use 
      the techniques of homological perturbation theory. Although the
treatment of BFV is basis dependent (individual constraints) and nominally
finite dimensional, I attempt to work more invariantly in terms of the ideal
generated by the constraints and take care to avoid assumptions of finite
dimensionality.  Although originally invented in the context 
of quantization, both BRST 
cohomology as they described it and the BFV-generalization are 
mathematically interesting in the `classical' setting.  The present paper is 
concerned only with the clasical setting but in the full generality of a 
first class ideal, in contrast to the paper of Kostant and Sternberg \c {KS} 
whose main interest is in  quantization issues for the case of an 
equivariant moment map and hence do not deal with the 
BFV-generalization nor with homological perturbation methods.
     \section{  Reduction}
           We have presented a geometric picture of reduction as 
referring to  $W \hookleftarrow V \rightarrow V/\cal F$.  There are a variety 
(pun intended) of difficulties with this approach.  The constraint locus 
$V$  fail to be a submanifold. Even if it is a
submanifold, the quotient $\hat V := V/\cal F$  may not be a manifold, in fact, 
may not even be Hausdorff.  (An intermediate situation of considerable
interest occurs with the quotient  $V/\cal F$  
being a stratified symplectic space \cite{LS}.)
	
When $(W,\omega)$ is a symplectic manifold with a smooth coisotropic 
submanifold, one of the nicest cases is called `regular', namely when the 
quotient $V/\cal F$ is a manifold and the projection $ V\to \hat V$ is a 
submersion. This implies further that
 $\omega|V$  has constant rank on  $TV$ (so that  $\omega|V$  is a presymplectic form on  $V$), and  $\cal F$ is an  involutive distribution  
given by $ker\  \omega|V$  which  is fibrating.
Then a standard argument, due essentially to E. Cartan \c{MW} or 
\c{GSte, Thm. 25.2}, 
shows that there exists a unique symplectic form $\hat \omega$ on  
$\hat V$  satisfying $\pi^* \hat \omega = \omega |V$.  The reduction 
of  $(W,\omega)$  is then the symplectic manifold  $(\hat V,\hat \omega)$
and the corresponding reduced Poisson algebra is  $C^{\infty}(\hat V)$
with the Poisson  bracket that is associated to $\hat \omega$.
           
In the ``singular'' case, when these conditions  
fail to hold, reduction in the above sense will not be 
      well defined.  Various definitions of reduction are possible, 
      depending upon which aspects of the theory are considered 
      primary.  (Of course, each such definition should agree with 
      regular reduction when both apply.)  Below we present two such 
      definitions (following \c {AGJ}), although there are undoubtedly 
      others. 

           The first type of reduction we shall consider is based upon 
      the notion of an ``observable''.  Following Bergman, we call a 
      function on  $W$  an {\bf observable} iff its Poisson bracket with 
each first class constraint is again a constraint, i.e., $h \in C^{\infty}(W)$
is an observable if and only if  $\{h,I\}\subset  I$.  Bergman emphasized 
      observables (rather than the points in  $V$  which are {\bf states}) 
      because observables represent measurable quantities.  (The 
      condition  $\{h,I\}\equiv 0$  on $V$  is a gauge invariance condition.)  
      The set   $\cal O (V)$  of observables forms a subalgebra of the 
      associative algebra  $C^{\infty} (W)$.
          
`` Dirac reduction'' takes two states  $x,y \in V$  to be 
physically equivalent iff they cannot be distinguished by observables.  
This amounts to defining  
an equivalence relation $\sim$  on  $V$  by  $x \sim y$
iff  $h(x) = h(y)$  for all observables  $h$.  
 The corresponding reduced space is  $\hat V = V/\sim$. 
The observables after 
reduction are identified with the elements of $\cal O (V)$ which are
fixed under the adjoint action of $I$ (with respect to Poisson bracket).
Since we are dealing with first class constraints, these observables
inherit a Poisson bracket.
           
\vskip1ex\noindent   
\noindent{\bf Example}:  Zero angular momentum in two dimensions.  
\vskip1ex\noindent   
Here  $W = T^* {\Bbb R}^2 \approx {\Bbb R}^2 \times {\Bbb R}^2 = \{(q,p)\}$
and the angular momentum is  $q \times  p = q_1p_2 - q_2p_1$  with constraint 
set  $V = \{(q,p)|q_1p_2 - q_2p_1 = 0\}$.  
      The foliation  $\cal F$  is in fact given by the orbits of
the standard circle action on 
${\Bbb R}^2$ lifted to  $T^*{\Bbb R}^2$.  The Dirac reduction can be identified 
with the symplectic orbifold ${\Bbb C}/Z_2$.
\vskip1ex
           
Sniatycki and Weinstein  \c{SW}  have defined an algebraic 
      reduction in the context of group actions and momentum maps which 
      is guaranteed to produce a reduced Poisson algebra but not 
      necessarily a reduced space of states (cf. \cite{W2}). (In contrast,
Kostant and Sternberg use the Marsden-Weinstein reduction \c{MW}.)
  The S-W (Sniatycki and Weinstein) reduced Poisson algebra 
is $ (C^\infty (W)/I)^G$   where  $V =  J^{-1}(0)$  
for some equivariant Poisson map  $J : W \rightarrow \frak g^*$  
(called a {\bf moment} map),
 equivariant with respect to a given $G$-action on $W,\quad \frak g$  being the 
Lie algebra of $G$.  (If  $G$  is compact and connected,
 $(C^{\infty} (W)/I)^G$   is 
isomorphic to the Dirac reduction  $C^{\infty} (W)^G /I^G$ .)  With hindsight, 
the generalization of S-W reduction to a general first class constraint ideal  
$I$ is obvious.  The issue of its suitability is not one of geometry 
      necessarily, but rather one of physics.  

The present paper grew out of the realization that the BFV construction could 
be regarded as a homological model which in degree zero models  the 
      I-invariants of  $C^{\infty} (W)/I$.  The whole construction turned out
in many cases  to
      be a model for the complex of longitudinal forms  $\Omega^*(V,\cal F)$.
From an algebraic geometric point of view, it is indeed natural to 
      define the observables on  $V$  by restriction of observables on  
      $W$, that is, to consider the quotient algebra  $C^{\infty} (W)/I$, which 
      corresponds to the algebra of smooth (in-the-sense-of Whitney) 
      functions on  $V$.  In physics, this is expressed by saying two 
  functions on  $W$  are weakly equal  $(f \approx  g)$  if their difference 
      vanishes on  $V$.
	
Now let us recast the problem in purely algebraic terms. Consider an
 arbitrary Poisson algebra $P$ with an ideal $I$ which is closed under the 
Poisson bracket.  Reduction is then achieved by passing to the $I$-invariant 
      subalgebra of  $P/I$.  Note that a class  $[g]$  is $I$-invariant if  
      $\{I,g\} \subset I$, equivalently, if  $\{\phi,g\} \approx   0$  for all 
constraints  $\phi \in I$. This subalgebra inherits a Poisson bracket even 
though  $P/I$  does not:  For  $f,g \in P$  and  
$\phi \in I,$ we have $ \{f + \phi,g\} = \{f,g\} + \{\phi,g\}$  
      where  $\{\phi,g\}$   need not belong to  $I$, but will if the class of
$g$  is $I$-invariant.

The Poisson algebra of invariants amounts to
 the quotient $N_P(I)/I$ where $N_P(I)$ denotes the normalizer of $I$ in $P$
 in the sense of Lie algebras; the ideal $I$ is a Poisson ideal in $N_P(I)$.

In this context, the analog of longitudinal forms are the alternating
multilinear-over-$P/I$ functions $h:I/I^2\otimes \dots\otimes I/I^2 \to P/I$ 
which again
form a graded commutative algebra, which we denote
$$ Alt_{P/I}(I/I^2, P/I).$$  We use $I/I^2$ because the corresponding 
Hamiltonian vector fields are restricted to $V$ in providing the 
foliation $\Cal F.$

The fact that  $I$  is a sub-Lie algebra of  $P$  but is not a
  Lie algebra  over  $P$  (the bracket is ${\Bbb R}$-linear but not $P$-linear)
      is a significant subtlety.  One way to handle this is to observe that 
$I/I^2$ inherits the structure of what 
Rinehart called an $(\Bbb R,P/I)$-{\bf  Lie algebra}.  This 
corresponds to what Herz \c{Hz} called a quasi-Lie algebra and what
Palais \c P called a $d$-Lie ring.
Since it is Rinehart's paper that establishes 
the relation to the geometry and was his major contribution in a 
tragically short career, we prefer to refer to the  Lie-Rinehart pair 
$(I/I^2, P/I)$.

\begin{df} \c R,\c P  A Lie-Rinehart pair $(L,A)$ over a ground ring $k$  
consists of a commutative $k$-algebra $A$ and
a Lie ring $L$ over $k$ which is a module over $A$ together 
with an $A$-morphism $\rho: L\to Der \ A$ such that
$$
[\phi,f\psi ] = (\rho (\phi)f)\psi + f[\phi,\psi ] \text{\ for\ }
\phi,\psi \in L, f \in A.
$$
\end{df}
Notice this is the condition satisfied by $L= I/I^2$ and $A=P/I$ with 
$\rho(\phi) f = \{\phi, f\}.$
      Hence we can consider the {\bf Rinehart complex}
$  Alt_{P/I}(I/I^2, P/I) $ with differential $d$ given by (3.1)
$$
            (dh)(\phi_0,...,\phi_q) = \sum_{i < j} (-1)^{i+j}
h([\phi_i,\phi_j ],...,\hat \phi_i,...,\hat \phi_j,...) + \sum_i 
(-1)^i \rho (\phi_i)h(...,\hat \phi_i,...).
$$

Realizing that $d$ is a derivation with respect to the usual product of
alternating functions, it is sufficient to know the above definition
for $q = 0 \text{\ and } 1.$ This differential given by
Rinehart \cite R is an obvious generalization of that of
Cartan-Chevalley-Eilenberg. 
      When $P/I$ is replaced by $P = C^{\infty}(W)$  and  $I/I^2$ by
the Lie algebra corresponding to  vector fields 
on $W$, the Rinehart complex becomes the de Rham complex of $W$. 
 As remarked by Stephen Halperin, 
the Rinehart complex  $Alt_{P/I} (I/I^2,P/I)$ is,
when $P=C^{\infty}(W)$ and $I$ is a first class ideal,  the complex
 $\Omega^* (V,\cal F)$  of longitudinal forms.
(See Huebschmann  \c {Hu1, Hu2, Hu5} for further applications of 
Rinehart's complex to 
Poisson algebras.)
           
This is the complex we wish to ``model''.  We will do this 
      using just the Poisson algebra structure of  $P$  and the sub-Lie 
      algebra and $P$-ideal  $I$, in contrast to the treatments of  \c{FHST}
      and  \c{D-V}  which retain some of the local manifold properties of  $W$.
     \section{  Differential Graded Commutative Algebras}
           One of the hallmarks of homological algebra is the use of 
      resolutions; for differential homological algebra, ``models'', in 
      the sense to be described, are more useful for many purposes.  
      For our approach to constrained Hamiltonian systems, one of the 
      basic objects is the deRham complex  $(\Omega^*(M),d)$  of differential 
      forms on a smooth manifold regarded as a DGCA (differential 
      graded commutative algebra):
             
$\Omega^*(M) = \{\Omega^p(M)\}$  where $\Omega^p(M)$  denotes 
the (real) vector space of differential p-forms, 

the wedge product  
$\omega \wedge \eta$ of forms gives $\Omega^*(M)$  the structure of a 
graded commutative algebra (over $\Bbb R$) : $\Omega^p \wedge \Omega^q \subset
\Omega^{p+q}$ with  $\omega \wedge \eta = (-1)^{pq}\eta \wedge \omega$, 

the exterior derivative  $d : \Omega^p \rightarrow \Omega^{p+1}$ is a graded 
derivation: $d(\omega \wedge \eta ) = d \omega \wedge \eta + (-1)^p  \omega
\wedge  d \eta$ and  $d^2  = 0$.
           
Another DGCA that plays an important role in mathematical 
      physics is the Cartan-Chevalley-Eilenberg complex  $(\Lambda \frak g^*,d)$
for the cohomology of a Lie algebra  $\frak g$.
Here, if  $\frak g$ is finite dimensional,  $\Lambda \frak g^*$
is usually interpreted as the exterior algebra  $E(\frak g^*)$  
on the $\Bbb R$-dual of 
$\frak g$, but, in general,  $\Lambda \frak g^*$   should be interpreted as
$Alt_{\Bbb R}(\frak g,{\Bbb R} )$, the algebra of alternating 
multilinear functions on  $\frak g$.  The coboundary  $d$  is given by 
(3.1) with  $\phi_i \in \frak g$.
           
Rational homotopy theory is much simpler than ordinary 
      homotopy theory because, for a large class of spaces, it is 
      completely equivalent to the homotopy theory of DGCAs over the 
      rationals  \cite{Q}.  Moreover, computations as well as 
      theoretical analysis can be carried out effectively in terms of 
      the models of Sullivan \cite{Su}.
\begin{df}  In the category of DGCAs over any $k$-algebra $P$, a 
  {\bf model} of a DGCA  $(A,d)$  is a morphism 
$\pi: ({\cal A},\partial) \rightarrow (A,d)$  of DGCAs such that  $\cal A$  
is free as a graded commutative algebra over $P$  
and $\pi^*: H({\cal A})\approx  H(A)$.
\end{df}
Here, 
%when dealing with algebras of finite type over $P$,
 free as a graded commutative algebra over $P$ means ${\cal A}$
is of the  form $P\otimes E(Z^{odd}) \otimes
k[Z^{even}]$ where $E$ = exterior algebra and $k[\ \ ]$ = polynomial
algebra and $Z$ is some free graded vector space of finite type.
Following the tradition in rational homotopy theory, the free graded
commutative algebra on a graded vector space $Z$
will be denoted $\Lambda Z.$
%More generally, $Alt_P(Z, \Lambda Y)$ is
% free as a graded commutative algebra over $P$. 
           
A major point of the Cartan-Chevalley-Eilenberg construction 
      in the case of a compact Lie group  $G$  is a natural map
      $(\Lambda \frak g^* ,d) \rightarrow \Omega^* (G)$  inducing a homology
isomorphism, i.e., a model for  $\Omega^* (G)$.
           
The main thrust of this paper is the construction of a differential
graded Poisson algebra which is, in many cases, a model 
      for the forms along the leaves of the constraint variety of a 
      first class system, in particular, $H^0$ will be isomorphic to 
the algebra of observables in the reduced sense: $(P/I)^I$ .

  \section{ Models for  $P/I$  and  $Alt_{P/I} (I/I^2,P /I)$  and the 
BRST generator}
Now we reverse the 
      procedure of  BFV  and first provide a model for  $P/I$  as a 
     $P$-module.  This model is a DGCA  $(P \otimes \Lambda\Psi,\delta )$ where
$\Psi$ is a graded vector space (in fact, negatively graded) and     
$\Lambda$ continues to denote the free graded commutative algebra.  
(This grading is the opposite of the usual convention in homological
algebra, but is chosen to correspond to the (anti-) ghost grading in the
physics literature and because we are modelling a DGCA of differential forms.)
 This model is constructed as follows:  Let $\Phi$ be the space of 
      $P$-indecomposables of  $I$ , i.e., $\Phi = I/\bar P I$  where  $\bar P$
is a complement to the constants ${\Bbb R} \subset  P$.  Let  $s\Phi$ denote 
a copy of $\Phi$  but regarded as having degree -1.  Let $\delta$ be the 
derivation of  $P \otimes \Lambda s \Phi$ determined by choosing a splitting
$\Phi \rightarrow I$ and factoring it as 
$\delta s: \Phi \rightarrow s\Phi \rightarrow I.$
(In terms of representatives  $\rho \in \Phi, \quad \delta\rho$ is  
$s^{-1}\rho$.)  In other words, $P \otimes \Lambda s\Phi$
is the Koszul complex for the ideal  $I$  in the 
commutative algebra  $P$  \c {Ko}, \c {Bo}.  If  $I$ is what is now known as a 
regular (at one time: Borel) ideal (an algebraic condition, but implied by
$I$   being  the defining ideal in  $C^{\infty}(W)$  for  $V = J^{-1}(0)$  
when  $0$  is a regular value of  $J : W \rightarrow {\Bbb R}^N )$, the Koszul
complex  $(P \otimes \Lambda s\Phi,\delta )$  is a model for  $P/I$.  For 
more general ideals, this fails, i.e., $H^i (P \otimes \Lambda s\Phi ,\delta)
\neq  0$  for some  $i \neq  0$.  The Tate resolution \c T kills this homology 
by systematically enlarging  $s\Phi$   to a graded vector space $\Psi$
and gives a model  $(P \otimes \Lambda\Psi,\delta)$  as desired.  We 
      refer to this model as  $K_I$ for brevity.  It is graded by the 
      grading on $\Psi$ extended multiplicatively, $\delta$   being still of 
      degree 1.
           
Now we wish to replace  $P/I$  by  $K_I$ in  \alt  with 
      the Rinehart generalization of the Cartan-Chevalley-Eilenberg 
      differential  $d$  and further alter it to a model which is itself a 
      (graded) Poisson algebra. The construction can be carried out quite
generally, but we succeed in showing it is a model in our sense
most easily in the case of a regular ideal, which obtains
 under reasonable geometric conditions.  Following the major theorem,
we describe a few other cases in which the model property also holds.
(Lars Kjeseth is continuing the purely algebraic version of this 
class of examples.  Kimura \c {Ki} has shown that for constraints which are not 
first class, the corresponding complex is NOT in general a model for the 
complex of forms along the leaves.)
\begin{th}
If $I$ is a first class ideal, there is a structure of differential
 graded Poisson algebra on $(\Lambda\Psi)^* \otimes  P \otimes \Lambda \Psi$
and a map of differential graded  Poisson algebras  
$$
   \pi  : ((\Lambda\Psi)^* \otimes  P \otimes \Lambda \Psi),\partial)
\rightarrow Alt_{P/I}(I/I^2, P/I)
$$
which induces an isomorphism on cohomology in degree zero.

      Here $\partial$ is $\delta + d +$ ``terms of higher order''  
in a sense to be made precise below.
\end{th}

The algebra structure on $ (\Lambda\Psi)^* \otimes  P \otimes \Lambda \Psi$
is that of the algebra of graded symmetric multilinear functions. 
The map $\pi$ is fairly straightforward.  Map  $P \otimes \Lambda\Psi
\rightarrow P/I$  by projection onto  $P$  and then by the quotient onto  $P/I$.
      Similarly project  $(\Lambda \Psi  )^*$  onto $(\Lambda s\Phi )^*$  
(recall $s\Phi$ is a summand of  $\Psi$ )  and then, identifying  
$(\Lambda  s\Phi )^*$
with  $Alt(\Phi ,R)$, map this to  $Alt_P (I,{\Bbb R} )$  by pulling back over 
the quotient  $I \rightarrow I/\bar PI = \Phi$. Finally, note the isomorphism 
of algebras  $ Alt_P(I, P/I) \iso$ \alt.  

We will construct the differential $\partial$ without any assumption on the 
ideal $I$ other than that it is first class.  The entire
construction $ ((\Lambda\Psi)^* \otimes  P \otimes \Lambda \Psi),\partial)$
we will denote by $X$.  In the  full generality of a first call ideal, 
we will show $H^i(X) =0$  for  $i < 0$
      and  $H^0(X)\approx (P/I)^I$ and moreover the isomorphism is given by the
  inclusion $(P/I)^I \hookrightarrow  P/I \hookrightarrow P \otimes \Lambda\Psi$
via the chosen splitting $P/I \hookrightarrow P$.  This then gives a
``no-ghost theorem'': $H^0(X)$  is represented completely by elements of
$P$  without any ghost (or  antighost) contributions from $\Lambda\Psi^*$ (or
$\Lambda\Psi$).  

For  $i > 0$, $H^i(X) $
      must be represented with ghosts. When this involves
only ghosts corresponding
      directly to constraints (i.e., elements of  $(s\Phi)^*$) 
      but no ghosts-of-ghosts,   ``geometrically'' we
      are looking at longitudinal forms.  It is only from the
      transverse (``gauge-fixed'') point of view that the ghosts inherit
      their name.

The key to the main theorem comes from the Hamiltonian and BRST formalisms.  
Let  $(\Lambda\Psi  )^* \otimes  P \otimes \Lambda\Psi$
be given a bigrading  $(r,s)$. Assuming  $P$  ungraded (see
\S 6 for the graded or super case), $P \otimes \Lambda\Psi$  is already 
(negatively) graded and this grading is  $s$, called the {\bf resolution} 
degree.  Then  $(\Lambda\Psi)^*$ inherits the dual (positive) grading  $r$, 
called the {\bf ghost} degree, adopting the term from the physics literature 
      (where the negative of the resolution degree is called the 
      anti-ghost degree).  The total degree is the sum  $r+s$  of the 
      ghost degree and the resolution degrees.  
Batalin, Fradkin, and Vilkovisky make  $X$  into a Poisson algebra by 
      extending the Poisson bracket on  $P$  to one on  $X$  by defining 
$$
\{h,\psi\} = h(\psi)\quad  {\rm for}\quad  h \in \Psi^ *,\quad \psi \in \Psi ,
$$
      all other brackets not determined by the derivation property 
      being set equal to zero.  This extended bracket is of total 
      degree zero, but mixed bidegrees. 
        \subsection{The BRST generator}
 The sought-for differential  $\partial$  is constructed to be of the 
form  $\partial = \{Q,\quad \}$  where  $Q$  is a formal sum of
terms  $Q_n$   defined by induction (on $n$).  
In physics, $Q$  is referred to as 
a BRST generator or operator,  in keeping with the philosophy mentioned in \S 2 
      with particular emphasis on the facts that  1) $\partial^2  = 0$  or 
      equivalently, $\{Q,Q\} = 0$ and  2) $Q$ contains a piece corresponding 
      to the Cartan-Chevalley-Eilenberg differential.
           
The proof of the existence of  $Q$  can be handled effectively 
by the ``step-by-step obstruction"  methods of homological perturbation theory  \c {G,GM,GSta,GLS,Hu6-9,HK}.  We adapt the details to this case, 
rather than appealing to the general theory.  
%(An alternate construction using the technique of Barnes 
%and Lamb\c{BL} has been proposed by Tatar and Tatar \c {TT}.)
We make crucial use of the filtration of  $X$  by the 
form or monomial degree, i.e., $(\Lambda^i\Psi)^* \otimes P \otimes \Lambda\Psi$
is the part of  $X$  of form degree  $i$, or equivalently, ``form degree $i$'' 
refers to an $i$-multilinear graded symmetric function from  
$\Psi$  to   $P \otimes \Lambda\Psi$.  The 
filtration is defined by: ${\cal F}^n  = {\cal F}^nX$  is the space of forms of 
degree  $> n$.  We use the strict inequality so that this filtration is 
multiplicative with respect to both parts of the Poisson algebra structure:
$$
{\cal F}^p {\cal F}^q \subset {\cal F}^{p+q+1} \subset {\cal F}^{p+q}      
\text{\ and\ }   \{{\cal F}^p,{\cal F}^q\} \subset {\cal F}^{p+q} .
$$
\indent  Start with  $Q_0: \Psi \rightarrow P \otimes \Lambda\Psi$ as 
the Koszul-Tate differential $\delta  $ restricted to $\Psi$.  As an element 
of  $X$, this  $Q_0$ is of total degree 1 and form degree 1, but  
$\{Q_0,\quad\}$  is a sum of two pieces, of form degree $0$  (namely 
$1 \otimes\delta$)  and of form degree 1.  Since the bracket restricts to 
the pairing (by evaluation) of  $(\Lambda\Psi)^*$  and  $\Lambda\Psi$, 
the term of form degree 1 includes the adjoint of  $\delta$   taking  
$Hom_P (\Lambda\Psi \otimes P,P)$  to itself.  The remainder of  $\{Q_0,\quad\}$
      is given by the original bracket (in $P$) of the coefficients of  
    $Q_0$ with elements of  $P$.
           
Since all our objects are at least vector spaces, the model 
      property of  $P \otimes \Lambda\Psi$ can be evidenced by a 
``contracting homotopy''  $s : P \otimes \Lambda\Psi \rightarrow
 P \otimes \Lambda\Psi $ of degree -1  such that  $s\delta + \delta s = 
1 - \bar\pi$   where  $\bar\pi : P \otimes  \Lambda\Psi \rightarrow P 
\rightarrow P/I \hookrightarrow  P \otimes \Lambda\Psi$ is given by  $\pi$    
composed with an  $\Bbb R$-linear splitting  $P/I \hookrightarrow P$. 
           
For any element  $R \in X$, let  $R^2$  denote  $1/2\{R,R\}$.  Now 
construct $ R_n  =  \sum_{i\leq n}  Q_i$   by induction so that
$$
\{R_n,R_n\} \in {\cal F}^{n+2} \quad \text { and} \quad
	\delta\{R_n,R_n\} \in {\cal F}^{n+3}.
$$
      Define  $Q_{n+1} = -s/2\{R_n,R_n\} = -sR_n^2$.

The following slightly complicated computation shows  $R_{n+1}$     
      satisfies the inductive assumption.
           
Both $\delta$ and  $s$  preserve the filtration, and from the way  
$Q_0$ is defined, $\{Q_0,\quad\}-1 \otimes\delta$    increases filtration.  
Start with
$$
R_{n+1}^2  =  (R_n  + Q_{n+1}^2)^2 =  R_n^2   - \{R_n,sR_n^2\} + (sR_n^2)^2 .
$$
      The last term  $(sR_n^2)^2 \in {\cal F}^{2n+4}$ since  
$sR_n^2 \in {\cal F}^{n+2}$ and  $2n+4 \geq n+4$.   On the other hand,
$$
 \{R_n,sR_n^2\} \equiv   (1 \otimes\delta)(sR_n^2 ) \quad mod \quad 
{\cal F}^{n+3},
$$
since  $R_n  = Q_0 + Q_1 + \dots $ and  the $\{Q_i, \quad \}$  for  $i > 0$  
increase filtration.  Thus
$$
\{R_n,sR_n^2\}\equiv -(1 \otimes s \delta)R_n^2 + R_n^2 \quad mod \quad 
{\cal F}^{n+3},
$$
      so
\begin{align} R_{n+1}^2 &\equiv -(1 \otimes s\delta) R_n^2 + R_n^2 \quad mod \quad
{\cal F}^{n+3}\\ 
        & \equiv 0 \quad mod \quad {\cal F}^{n+3}
\end{align}
by the assumption on $\delta  R_n^2$  .

Similarly 

\begin{align} \delta R_{n+1}^2 &\equiv \delta R_n^2 - \delta\{R_n, sR_n^2\} + \delta (sR_n^2)^2 \\ 
	& \equiv \delta R_n^2 \quad mod \quad {\cal F}^{n+4}.
\end{align} 

      Now we need to commute $\delta$ with  $\{R_n,\quad\}$.  Since  
$\{R_n,\quad \} - 1 \otimes \delta$ increases filtration by at least one, 
its square does so by at least two.  Thus
$$
\{R_n,\{R_n,\quad \}\} - \{R_n, 1 \otimes\delta\} - 1\otimes\delta \{R_n,\quad\}
$$
applied to  $sR_n^2$ is of filtration at least $n+4$.  Now the graded Jacobi 
identity gives
$$
                         2\{R_n,\{R_n,\quad\}\} = \{\{R_n,R_n\},\quad \}
$$
      which increases filtration by  $n+2$, thus
\begin{align} \delta R_{n+1}^2 &\equiv \delta R_n^2  + \{R_n,\delta sR_n^2\}\quad mod\quad {\cal F}^{n+4}\\ 
         & \equiv \delta R_n^2   + \{R_n ,R_n^2\} - \{R_n ,s\delta R_n^2\} 
\quad mod \quad {\cal F}^{n+4}\\ 
    & \equiv \delta R_n^2 - \delta R_n^2 \quad mod \quad {\cal F}^{n+4}
\end{align}
since  $\{R_n ,s \delta R_n^2\}  \equiv \delta s\delta R_n^2 \quad mod 
{\cal F}^{n+4}$  and  $\{x,\{x,x\}\} = 0$  for  $x$  of any degree (over a 
field of characteristic not equal to 3).

Thus we have constructed a differential graded Poisson algebra for any 
coisotropic ideal.  Where possible, we will show that we have a model for
$Alt_{P/I}(I/I^2, P/I)$ by the usual 
      techniques of comparison in homological perturbation theory, 
      namely comparison of spectral sequences.  In one final case,
we can do this locally but appeal to a geometric arguement to patch
the local results.
After establishing that, we will look at issues involving choices (possibly 
non-minimal) choices of generators (constraints) for the ideal $I$.

           From the definition of the filtration  ${\cal F}^p$, the associated
      graded  $E_0 (X)$  is isomorphic to  $(\Lambda\Psi)^*  \otimes  P \otimes 
\Lambda\Psi$.  To analyze  $d_0$, notice that since $s$ preserves the form 
degree, $Q_{i+1}= -sR_i^2 \in \cal F^{i+2}$ and hence $\{Q_i,\ \}$ increases 
filtration by at least $i$.  As mentionned earlier, $\{Q_o,\ \} - 1\otimes 
\delta$ also increases filtration so $d_0$ is (induced by) the Koszul 
differential $\delta$.  Thus 
$$
   E_1 (X) \approx (\Lambda \Psi)^* \otimes P/I \approx Alt_{\Bbb R}(\Psi, P/I),
$$
and $E_1(X)$ is concentrated in anti-ghost degree 0; the spectral sequence
necessarily collapses from $E_2$.  To determine 
$H^0(X) \approx E_2^{0,0}$, we need only analyze $d_1$ on $\Phi$. For 
$h\in P/I,$ consider $\{Q_0,h\}: I \to P/I$.  It is given by $\{Q_0,h\}(\phi) 
= \{\phi, h\}$ for $\phi\in I.$  Thus $H^0(X)$ is isomorphic to the $I$-
invariants of $P/I$.

When the ideal $I$ is regular, $\Psi = s\Phi$ and we can analyze $d_1$ on 
$\Lambda s\Phi$ similarly.  
For example, for
$h:I \rightarrow P/I$, consider  
$\{Q_0,h\}:I \wedge I \rightarrow P/I$.  It is given by  
$\{Q_0,h\}(\phi_1,\phi_2) = \{\phi_1,h(\phi_2)\} - \{\phi_2,h(\phi_1)\}$  while 
$\{Q_1,h\}(\phi_1,\phi_2) = - 1/2\{s\{Q_0,Q_0\},h\}(\phi_1,\phi_2) =
-h\{\phi_1,\phi_2\}$.  (At this point, one appreciates the facility of 
non-invariant description in terms of a generating set of constraints  
$\{\phi_{\alpha}\}$ for $I$ and a dual set $\{\eta^{\alpha}:I \rightarrow P\}$.)
           
Thus we see  $d_1$ (up to sign) looks like the Rinehart 
      generalization of the Cartan-Chevalley-Eilenberg differential.  
It is this identification of  $(E_1,d_1)$  which motivates the name   
       {\bf BRST generator} for $Q$.

Now to make the comparison with the complex of longitudinal forms,
    since $\Phi$ is defined as a quotient of $I$, there is the induced 
chain map

$$
 \pi:X \to  Alt_k(\Psi, P/I)\to Alt_k(s\Phi, P/I) \to Alt_P (I,P/I) \approx Alt_{P/I} (I/I^2,P /I)
$$
as described above.  In the regular case all maps except $\pi$ are isomorphisms.
For the constrained Hamiltonian setting with which we began,
      in which  $P$  is  $C^{\infty}(W)$, we have  identified \alt
with the longitudinal forms of the foliation  ${\cal F}$
  of  $V$  and  $d_1$ with the exterior derivative ``along the leaves''.
\begin{th}
If $I$ is a regular first class ideal in $C^{\infty}(W)$, the map
$\pi$ induces an isomorphism $H(X) \approx H(\Omega(V, \cal F))$.
\end{th}

When $I$ is not regular, we still have the map but in general lack
sufficient information to conclude an isomorphism in cohomology.

Now the physicists do not work with the ideal explicitly but rather with 
a set of constraints, which is a set (not necessarily minimal) of 
generators for the ideal.  The corresponding BFV construction starts with
$\Phi$ as the vector space spanned by the constraints, rather than
with $I/{\overline P} I$.  In certain cases, even though the constraints
do not form a regular sequence, we can still make the
identification of $H(X)$ with  $H(\Omega(V, \cal F))$.

The redundant case:  The set of constraints may be reducible in a trivial 
way; a proper subset may consist of  a regular sequence  of generators.  
Then we can split $\Phi$ as $\Xi\oplus \Upsilon$ where $\Xi$ is 
the span of the minimal subset and $\Upsilon$ is spanned by the complementary 
subset.  The Koszul-Tate resolution of $P/I$ splits as the Koszul resolution 
determined by $\Xi$ tensored with a contractible DCGA.  Then 
$Alt(\Psi,\ )$ splits similarly and the BRST generator can be constructed 
first in the $\Xi$ part and then extended so the results will be the same 
as when using $\Phi= I/{\overline P} I$.

In particular, if the constraints are given by an equivariant 
moment map $J : W \rightarrow \frak g^*$ where G acts by 
symplectomorphisms but with kernel H, then 
$I/{\overline P} I$ is isomorphic to $\frak g/\frak h$ but the span  
of the constraints would be isomorphic to $\frak g$. Here choose a 
splitting $\Xi\oplus \Upsilon$ such that $\Upsilon = \frak h$ and
$\Xi \approx \frak g/\frak h,$
then proceed as in the redundant case.
           
In \c {FHST} and \c{HT},  the setting is specifically that of
a symplectic manifold (phase space) with a constraint submanifold
(``surface'') and moreover the assumption is made that locally
the constraints can be separated into ``independent constraint functions''
and dependent ones which can be expressed as functional linear
combinations of the independent ones with coefficients which are 
regular in a neighborhood of the constraint submanifold.  Thus locally
we are in the redundant case so identities involving the globally
defined BRST generator and comparisons with the complex of forms along
the leaves can be verified locally; we again have
$H(X) \approx H(\Omega(V, \cal F))$.

Finally, the construction of  $\partial$  and of  $Q$  involves 
a choice of contracting homotopy  $s$  and implicitly of a choice of 
splitting  $P/I \hookrightarrow P$.  A change in  $s$  produces changes in  
$\partial$  but not in the homotopy type of  $(X,\partial)$  as DGCA.  
Moreover the change in  $s$  can be realized by an automorphism of $\Lambda\Psi$
and the induced one on $\Lambda\Psi^*$. This is an example of  what is known as 
a canonical transformation, a basic automorphism of any Hamiltonian system.
      \section{ Generalizations:  Infinite dimensions and super algebras}
           If  $I$ is regular and finitely generated over  $P$  (so 
$\Phi$ is finite dimensional over ${\Bbb R})$, $Alt_P(\Phi,P \otimes  s\Phi)$  
is finitely generated as a $P$-module and  $Q_n = 0$  for sufficiently large  
$n$.  If  $I$  is finitely generated but not regular, $\Psi$ may easily be 
infinite dimensional, though finite in each grading, and so all  $Q_n$   
may be non-zero.
           
More importantly, there are many examples occurring in physics 
      (field theory) in which $\Phi$ is itself infinite dimensional.  That 
      is why we have been careful to emphasize  $Alt$  or to take the 
      dual of $\Lambda\Psi$ rather than $\Lambda(\Psi^*)$.  Actually both 
physical and mathematical considerations (cf. Gelfand-Fuks cohomology) suggest 
      that the alternating functions might better be restricted to 
      being continuous in an appropriate topology.
           
Early in the development of Batalin, Fradkin and 
      Vilkovisky's approach, attention was called to the generalization 
      to a super-Poisson algebra  $P = P_0 \oplus P_1$   with super constraints.
      This means that  $P$ is a GCA (graded by  $Z/2 = \{0,1\}$) with a 
      graded bracket $\{\quad,\quad\}$:
\begin{align}  &P_0 \otimes P_0  \rightarrow P_0\\ 
	&P_0 \otimes P_1  \rightarrow P_1\\ 
	&P_1 \otimes P_1  \rightarrow P_0
\end{align}
      with graded anticommutativity, graded Jacobi identity, and graded 
      derivation property (Leibnitz rule):
$$                                             
                      \{f,gh\} = \{f,g\}h + (-1)^{|f||g|} g\{f,h\}
$$
 where  $f \in  P_{|f|}, \ g \in P_{|g|}$.
           
It has long been known in algebraic topology how to 
      generalize the construction of models such as the Koszul-Tate 
      complex or the Chevalley-Eilenberg complex to the graded setting, 
      e.g.,$\Psi$ is now a graded vector space and  $s\Phi$ is an isomorphic 
      copy of a $\Phi$ regraded down by $1$ so that $\delta$ is still of 
degree $1$.  The use of $\Lambda$ to denote the free graded commutative 
algebra on a 
      graded vector space means that the only necessary change in our 
      treatment is to specify the resolution degree as the one implied 
      by the degree on  $s\Phi$ with $\delta$ being of resolution degree   
$1$.  Notice this is not the same as ignoring the internal grading on  
      $s\Phi$ and just counting the algebraic degree.  (It is spelled out in 
      [GS1] for example.)  From there on, the signs take care 
      of themselves if we follow the usual conventions, introducing a 
      sign  $(-1)^{pq}$ whenever a term of total degree  $p$  is pushed past 
      a term of total degree  $q$.
\section{ References}

      [AGJ]      J.M. Arms, M. J. Gotay and G. Jennings, (Geometric and
                 Algebraic) Reduction for Singular Momentum Maps, Adv. in Math.
		79 (1990) 43-103.
\vskip.5ex
      [BF]       I.A. Batalin and E.S. Fradkin, A generalized canonical
                 formalism and quantization of reducible gauge theories,
                 Phys. Lett. 122B (1983) 157-164.
\vskip.5ex
      [BV1]      I.A. Batalin and G.S. Vilkovisky, Existence theorem for
                 gauge algebra, J. Math. Phys. 26 (1985) 172-184.
\vskip.5ex
      [BV2]      I.A. Batalin and G.S. Vilkovisky, Quantization of gauge
                 theories with linearly dependent generators, Phys. Rev. D
                 28 (1983) 2567-2582.
\vskip.5ex
      [BV3]      I.A. Batalin and G.S. Vilkovisky, Relativistic S-matrix
                 of dynamical systems with boson and fermion constraints,
                 Phys. Lett. 69B (1977) 309-312.
\vskip.5ex
      [BRS]	 C. Becchi and A. Rouet and R. Stora,  Renormalization of the 
		 abelian {H}iggs-{K}ibble model, Commun. Math. Phys. 42 (1975)
		 127-162.

%\vskip.5ex
      %[B]	    . Bergman, title, Phys. Rev. 98 (1955) 531-
\vskip.5ex
      [Bo]       A. Borel,  Sur la cohomologie des espaces fibr\'es principaux
et des espaces homogenes de groupes de Lie compacts, Annals of Math. 57 (1953)
115-207.
\vskip.5ex
      [BM]       A.D. Browning and D. McMullen, The Batalin, Fradkin, 
                 Vilkovisky formalism for higher order theories, J. Math.
                 Phys. 28 (1987) 438-444.
\vskip.5ex
      [CE]       C. Chevalley and S. Eilenberg, Cohomology theory of Lie
                 groups and Lie algebras, Trans. Amer. Math. Soc. 63 
                 (1948) 85-124.
\vskip.5ex

      [D]        P.A.M. Dirac, {\sl Lectures on Quantum Mechanics}, Belfer 
                 Graduate School Monograph Series 2 (1964).
\vskip.5ex
      [D-V]      M. Dubois-Violette, Systemes dynamiques constraints:
l'approche homologique, Ann. Inst. Fourier, Grenoble 37 (1987) 45-57.
\vskip.5ex
      [FF]       E.S. Fradkin and T.E. Fradkina, Quantization of relativistic 
system with boson and fermion first- and second-constraints, Phys. Lett.
72B (1978) 343-348.
\vskip.5ex
      [FHST]     J. Fisch, M. Henneaux, J. Stasheff and C. Teitelboim, 
                 Existence, uniqueness and cohomology of the classical 
                 BRST change with ghosts of ghosts, Comm. Math. Phys. 120 (1989)
379-407.
\vskip.5ex
      [FV]       E.S. Fradkin and G.S. Vilkovisky, Quantization of
                 relativistic systems with constraints, Phys. Lett. 55B
                 (1975) 224-226.
\vskip.5ex
      [GIMMSY]   M.J. Gotay, J. Isenberg, J.E. Marsden, R. Montgomery, J.
                 Sniatycki and P. Yasskin,  { \it Constraints and momentum
                 mappings in relativistic field theory}, preprint.
\vskip.5ex
      [G]       V.K.A.M. Gugenheim, On a perturbation theory for the 
                 homology of a loop space, J. Pure and Appl. Alg. 25 
                 (1982) 197-205.
\vskip.5ex
      [GLS]	V.K.A.M. Gugenheim, L. Lambe and J. Stasheff, 
		Algebraic aspects of Chen's twisting cochain, 
		Ill. J. M. 34  (1990) 485-502.
\vskip.5ex
      [GM]       V.K.A.M. Gugenheim and J.P. May, {\sl On the theory and
                 application of torsion products}, Memoirs Amer. Math. Soc.
                 142 (1974).
\vskip.5ex
      [GSta]      V.K.A.M. Gugenheim and J.D. Stasheff, On perturbations
                 and $A_{\infty}$-structures, Bull. Soc. Math. Belg. 38 (1986)
                 237-245.
\vskip.5ex
      [GSte]      V. Guillemin and S. Sternberg, {\sl Sympletic Techniques in
                 Physics}, Cambridge University Press, Cambridge 1984.
\vskip.5ex
      [HH]	 J. L. Heitsch and S. E. Hurder, Geometry of foliations,
		 J. Diff. Geom. 20 (1984) 291-309.
\vskip.5ex
      [H]        M. Henneaux, Hamiltonian form of the path integral for
                 theories with a gauge freedom, Phys. Rev. 126 (1985) 
                 1-66.
\vskip.5ex
      [HT]	M. Henneaux and Claudio  Teitelboim,{\it Quantization of 
		gauge systems}, Princeton University Press, Princeton, N.J. 1992. 
\vskip.5ex
      [Hu1]	 J. Huebschmann, Poisson cohomology and quantization, 
		 J. fuer die Reine und Angewandte Mathematik 408 (1990) 57-113
 
\vskip.5ex
[Hu2] On the quantization of Poisson algebras, 
 in {\it Symplectic Geometry and mathematical physics, actes du colloque en
 l'honneur de J.M. Souriau}, P. Donato, C. Duval, J. Elhadad, G. M.
 Tuynman, eds., Progress in Math. vol 99, Birkhaeuser, Basel, 
 Boston, Berlin, (1991) pp. 204--233 
\vskip.5ex
 
 [Hu3] Extensions of Lie Rinehart algebras. I. The Chern-Weil construction,
 Preprint 1989
\vskip.5ex
 
[Hu4] Extensions of Lie-Rinehart algebras. II. The spectral sequence and
 the Weil algebra, Preprint 1989
 
\vskip.5ex
[Hu5] Graded Lie-Rinehart algebras, graded Poisson algebras, and BRST-quantization. I. The finitely generated case.
 Preprint 1990.

\vskip.5ex
      [Hu6]	J. Huebschmann, Perturbation thoery and small methods 
		for the chains of certain induced fibre spaces, 
		Habilitationsschrift, Universit\"at Heidelberg  (1984)

\vskip.5ex
	[Hu7]	The cohomology of $F\Psi^q$. The additive structure,
	        J. Pure and Appl. Alg. 45 (1987) 73-91.

\vskip.5ex
	[Hu8]	Perturbation theory and free resolutions for nilpotent groups 
		of class 2, J. Alg. 126 (1989) 348-399
\vskip.5ex
	[Hu9]	Minimal free multi models for chain algebras, preprint  1990 

\vskip.5ex
      [HK]	 J. Huebschmann and T. Kadeishvili, Small models for 
		 chain algebras, Math. Z. 207 (1991) 245-280. 
\vskip.5ex
      [Hz] 	 J. Herz, Pseudo-alg\`ebres de Lie, C.R. Acad. Sci. Paris 236 
		 (1953) 1935-1937.
\vskip.5ex
      [Ki] 	 T. Kimura, Generalized classical BRST and reduction of 
		 Poisson manifolds, Comm. Math. Phys. 151(1993) 155-182
\vskip.5ex
      [KS]	 B. Kostant and S. Sternberg, Symplectic reduction, BRS 
		 cohomology and infinite-dimensional Clifford algebras, 
		 Annals of Physics 176 (1987) 49-113.
\vskip.5ex
      [Ko]        J.-L. Koszul, Sur un type d'algebres differentielles en 
		 rapport avec la transgression, {\it Colloque de Topologie}, 
		 Bruxelles (1950), CBRM, Liege
\vskip.5ex
      [LS]	 R. Sjamaar and E. Lerman, Stratified symplectic spaces 
		 and reduction, Annals of Math. 134 (1991) 375-422.
\vskip.5ex
      [MW]       J.E. Marsden and A. Weinstein, Reduction of symplectic
                 manifolds with symmetry, Rep. Math. Phys. 5 (1974)
                 121-130.
\vskip.5ex
       [P]	  R. Palais, The cohomology of Lie rings, Proc. Symp. Pure Math. III (1961)  130-137, AMS.	
\vskip.5ex
      [Q]        D. Quillen, Rational homotopy theory, Annals of Math. 90 
(1969) 205-295.
\vskip.5ex
      [R]        G. Rinehart, Differential forms for general commutative
                 algebras, Trans. Amer. Math. Soc. 108 (1963) 195-222.
\vskip.5ex
      [SW]       J. /'Sniatycki and A. Weinstein, Reduction and
                 quantization for singular momentum mappings, Lett. Math.
                 Phys. 7 (1983) 155-161.
\vskip.5ex
      [S1]       J.D. Stasheff, Constrained Hamiltonians: A homological
                 approach, Proc. Winter School on Geometry and Physics,
                 Suppl. Rendiconti del Circ. Mat. di Palermo, II, 16
                 (1987) 239-252.
\vskip.5ex
      [S2]       J.D. Stasheff, Constrained Poisson algebras and strong
                 homotopy representations, Bull. Amer. Math. Soc. (1988)
                 287-290.
\vskip.5ex
      [Su]       D. Sullivan,  Infinitesimal computations in topology, 
		 Publ.\ Math. IHES.\ vol.\ 47 (1977)  269-331.
\vskip.5ex
      [T]        J. Tate, Homology of Noetherian rings and local rings,
                 Ill. J. Math. 1 (1957) 14-27.
\vskip.5ex
      [Ty]	 I. V. Tyutin,  Gauge invariance in field theory and statistical
		 physics in operator formulation (in {R}ussian),  Lebedev Physics
		 Inst. preprint 39 (1975)

\vskip.5ex
      [W]	 A. Weinstein, Coisotropic calculus and Poisson groupoids,
		 J. Math. Soc. Japan 40 (1988) 705-727.
\vskip.5ex
      [W2]	 A. Weinstein, Symplectic V-manifolds, periodic orbits of 
		 Hamiltonian systems, and the volume of certain Riemannian 
		 manifolds, Comm. Pure Appl. Math 30 (1977) 265-271.
\noindent

%\bibliographystyle{amsplain}
%\bibliography{op}
%\end{document}

\begin{ack}
The author would like to thank the University of Pennsylvania
for hospitality during his leave  (and many summers) 
and Lehigh University for a subsequent
visiting appointment. 
\end{ack}
\end{document}